\documentclass[a4paper,12pt]{amsart}
\usepackage{amssymb}
\usepackage{ifthen}
\usepackage{graphicx}
\usepackage{float}
\usepackage{caption}
\usepackage{subcaption}
\usepackage{cite}
\usepackage{amsfonts}
\usepackage{amscd}
\usepackage{amsxtra}
\usepackage{mathrsfs}
\usepackage[usenames]{color}

\newtheorem{thm}{Theorem}
\newtheorem{cor}{Corollary}
\newtheorem{lem}{Lemma}

\newtheorem{rem}{Remark}
\newtheorem{example}{Example}
\newtheorem{defn}{Definition}
\newtheorem{prob}{Problem}

\newtheorem{conj}{Conjecture}
\theoremstyle{definition}

\newcounter {own}
\def\theown {\thesection  .\arabic{own}}

\newenvironment{pf}[1][]{%
 \vskip 3mm
 \noindent
 \ifthenelse{\equal{#1}{}}%
  {{\slshape Proof. }}%
  {{\slshape #1.} }%
 }%
{\qed\bigskip}

\newcounter{alphabet}
\newcounter{tmp}
\newenvironment{Thm}[1][]{\refstepcounter{alphabet}%
\bigskip%
\noindent%
{\bf Theorem \Alph{alphabet}}%
\ifthenelse{\equal{#1}{}}{}{ (#1)}%
{\bf .} \itshape}{\vskip 8pt}

\makeatletter
\newcommand{\Ref}[1]{\@ifundefined{r@#1}{}{\setcounter{tmp}{\ref{#1}}\Alph{tmp}}}
\makeatother

\newenvironment{Lem}[1][]{\refstepcounter{alphabet}%
\bigskip%
\noindent%
{\bf Lemma \Alph{alphabet}}%
{\bf .} \itshape}{\vskip 8pt}

\newcommand{\IR}{{\mathbb R}}

\newcommand{\IC}{{\mathbb C}}
\newcommand{\ID}{{\mathbb D}}

\newcommand{\K}{{\mathcal K}}

\def\be{\begin{equation}}
\def\ee{\end{equation}}

\newcommand{\bee}{\begin{enumerate}}
\newcommand{\eee}{\end{enumerate}}

\newcommand{\blem}{\begin{lem}}
\newcommand{\elem}{\end{lem}}
\newcommand{\bthm}{\begin{thm}}
\newcommand{\ethm}{\end{thm}}
\newcommand{\bcor}{\begin{cor}}
\newcommand{\ecor}{\end{cor}}
\newcommand{\beg}{\begin{example}}
\newcommand{\eeg}{\end{example}}
\newcommand{\begs}{\begin{examples}}
\newcommand{\eegs}{\end{examples}}
\newcommand{\bdefe}{\begin{defn}}
\newcommand{\edefe}{\end{defn}}
\newcommand{\bprob}{\begin{prob}}
\newcommand{\eprob}{\end{prob}}
\newcommand{\bei}{\begin{itemize}}
\newcommand{\eei}{\end{itemize}}

\newcommand{\bcon}{\begin{conj}}
\newcommand{\econ}{\end{conj}}
\newcommand{\bcons}{\begin{conjs}}
\newcommand{\econs}{\end{conjs}}
\newcommand{\bprop}{\begin{propo}}
\newcommand{\eprop}{\end{propo}}
\newcommand{\br}{\begin{rem}}
\newcommand{\er}{\end{rem}}
\newcommand{\brs}{\begin{rems}}
\newcommand{\ers}{\end{rems}}
\newcommand{\bo}{\begin{obser}}
\newcommand{\eo}{\end{obser}}
\newcommand{\bos}{\begin{obsers}}
\newcommand{\eos}{\end{obsers}}
\newcommand{\bpf}{\begin{pf}}
\newcommand{\epf}{\end{pf}}
\newcommand{\ba}{\begin{array}}
\newcommand{\ea}{\end{array}}
\newcommand{\beq}{\begin{eqnarray}}
\newcommand{\beqq}{\begin{eqnarray*}}
\newcommand{\eeq}{\end{eqnarray}}
\newcommand{\eeqq}{\end{eqnarray*}}

\newcommand{\ds}{\displaystyle}

\font\ff=eusm10
\def\K{\hbox{\ff K}}

\newcounter{minutes}
\divide\time by 60
\newcounter{hours}
\multiply\time by 60 \addtocounter{minutes}{-\time}

\begin{document}
\bibliographystyle{amsplain}
\title[Riesz-Fej\'er inequalities for harmonic functions]
{Riesz-Fej\'er inequalities for harmonic functions}

\thanks{
File:~\jobname .tex,
          printed: \number\year-\number\month-\number\day,
          \thehours.\ifnum\theminutes<10{0}\fi\theminutes}

\author{Ilgiz R Kayumov}
\address{I.R. Kayumov, Kazan Federal University, Kremlevskaya 18, 420 008, Kazan. Russia.}
\email{ikayumov@kpfu.ru}

\author{Saminathan Ponnusamy
}
\address{S. Ponnusamy
Indian Institute of Technology Madras, Department of Mathematics, Chennai –600 036, India. }
\email{samy@isichennai.res.in, samy@iitm.ac.in}

\author{Anbareeswaran Sairam Kaliraj}

\address{ A. Sairam Kaliraj,
Indian Institute of Technology Ropar, Nangal Road, Rupnagar, Punjab - 140001, India. }
\email{sairamkaliraj@gmail.com, sairam@iitrpr.ac.in}

%


\date{\today  
}

\begin{abstract}
In this article, we prove the Riesz - Fej\'er inequality for complex-valued harmonic functions in the harmonic Hardy space ${\bf h}^p$ for all $p > 1$. The result is sharp for $p \in (1,2]$. Moreover, we prove two variant forms of Riesz-Fej\'er inequality for harmonic functions, for the special case $p=2$.
\end{abstract}
\thanks{ }

\maketitle
\pagestyle{myheadings}
\markboth{I.R. Kayumov, S.Ponnusamy, and A. Sairam Kaliraj}{Riesz-Fej\'er type inequalities for harmonic functions}
\section{Introduction and Main Results}\label{KPSSec1}

Let $\ID$ be the open unit disk in the complex plane, i.e. ${\mathbb D}=\{z \in {\mathbb C}:\, |z|<1\}$ and let
${\mathcal A}$ denote the class of all analytic functions $f$ defined on $\ID$.  For $f \in \mathcal A$,
the integral means $M_p(r, f)$ is defined as
$$
M_p(r, f) = \left( \frac{1}{2\pi}\int_{0}^{2\pi}|f(r e^{i\theta})|^p d\theta \right)^{1/p}, ~~~0 < p < \infty .
$$
The classical Hardy space $H^p$, $0 < p < \infty $, consists of all analytic functions $f: \ID \to \IC$ such that $M_p(r, f)$ remains
bounded as $r \to 1^-$.
The study of $H^p$ spaces attracted the attention of many mathematicians, as it deals with the important problems in function theory such as the existence of
the radial limit in almost all directions, growth of the absolute value of functions, bounds on the coefficients and so on.  For more details on $H^p$ spaces, one can refer to the books of \cite{Duren:Hardy, Koosis}. One of the celebrated results on $H^p$ spaces by Riesz -- Fej\'er is the following inequality:

\begin{Thm}\label{theA}{\rm \cite[Theorem 3.13]{Duren:Hardy}}
If $f \in H^p$ $(0 < p < \infty)$, then the integral of $|f(x)|^p$ along the segment $-1 \leq x \leq 1$ converges, and
\be\label{KPSeq1}
\int_{-1}^{1}|f(x)|^p dx \leq \frac{1}{2}\int_{0}^{2\pi}|f(e^{i\theta})|^p d\theta ,
\ee
where $f(e^{i\theta})$ stands for the radial limit of $f$ on the unit circle. The constant $1/2$ is best possible.
\end{Thm}

The above theorem has a nice geometry: When the unit disk $\ID$ is mapped conformally onto the interior of a rectifiable Jordan curve $C$,
the image of any diameter has at most half the length of $C$.

Theorem \Ref{theA} has been generalized in several settings over the years. For example,
Beckenbach \cite{Beckenbach} proved that the inequality \eqref{KPSeq1}
remains valid if $|f|^p$ is replaced by a non-negative function whose logarithm is subharmonic. Validity of inequality \eqref{KPSeq1}
under weak regularity assumptions and generalizations of it may be seen from \cite{Beckenbach,Calderon,Huber} and the references therein. However, not much is known about
Riesz-Fej\'er inequality for complex-valued harmonic functions defined on the unit disk $\ID$. In this article, we prove Riesz-Fej\'er inequality for complex-valued harmonic functions
and the results are sharp for the cases $1 < p \leq 2$.

Let ${\mathcal H}$ denote the class of all complex-valued harmonic functions $f=u+iv$, where $u$
and $v$ are real-valued harmonic functions in $\ID$. It is easy to see that the function
$f$ has the unique decomposition $f=h+\overline{g}$ in the unit disk ${\mathbb D}$, where $h$ and $g$ are analytic
in $\mathbb{D}$ with the normalization $g(0)=0$. The normalization $g(0)=0$ do not affect any result pertaining
to the integral means as we can always induct $\overline{g(0)}$ into $h(0)$.
For $0<p<\infty$, denote by ${\bf h}^p$, the space of all complex-valued harmonic functions $f \in \mathcal{H}$
satisfying the condition $M_p(r, f)$ is uniformly bounded with respect to $r$. Each function $f \in {\bf h}^p$ has a radial limit almost everywhere. Throughout this article, we use  $f(e^{i\theta})$ to denote the radial limit of $f$ on the unit circle. For $f \in {\bf h}^p$, define $||f||_p$ as
$$||f||_p:=\left( \frac{1}{2\pi}\int_{0}^{2\pi}|f(e^{i\theta})|^p d\theta \right)^{1/p} < \infty.
$$
There are many results in the literature that  are more or less directly connected with
the Riesz-Fej\'er inequality \eqref{KPSeq1}. Riesz and Zygmund proved the harmonic analog of \eqref{KPSeq1}
for the special case $p=1$ which has the same geometric interpretation as that of their analytic counterparts.

\begin{Thm}[Riesz-Zygmund inequality; See also {\rm \cite[Theorem 6.1.7]{Pavbook}}]\label{the_RieZyg}
If $f \in {\mathcal H}$ such that $\partial f(r e^{i \theta})/ \partial \theta \in {\bf h}^1$, then
\be\label{KPSeq2}
\int_{-1}^{1}\left|\frac{\partial f(re^{it})}{\partial r}\right| dr \leq \frac{1}{2} \sup_{0 < r < 1} \int_{0}^{2\pi}\left|\frac{\partial f(re^{i\theta})}{\partial \theta}\right| d\theta .
\ee
The constant $1/2$ is sharp.
\end{Thm}

As a Corollary to Theorem \Ref{the_RieZyg}, one could easily get the following result: When $f$ is a harmonic diffeomorphism from the unit disk $\ID$
onto a Jordan domain with rectifiable boundary of length $C$, then the image of any diameter under $f$ has at most half the length of $C$.

In this article, we prove the sharp form of Riesz -- Fej\'er inequality for $f \in {\bf h}^p$ for $p \in (1,2]$ and it is as follows:

\bthm\label{KPSthm1} Let $p>1$. If $f \in {\bf h}^p$ then the following inequality holds:
$$\int_{-1}^{1}|f(x e^{i t})|^p dx \leq  A_p  \int_{0}^{2\pi}|f(e^{i\theta})|^p d\theta, ~~\mbox{ for all }~~ t \in \IR.
$$
where
$$A_p := \left\{
    \begin{array}{ll}
        \ds \frac{1}{2}\sec^{p} \left(\frac{\pi}{2 p}\right)  & \mbox{if } 1 < p \leq 2 \\
        \ds 1 & \mbox{if } p \geq 2.
    \end{array} \right.
$$ The estimate is sharp for all $p \in (1,2]$.

\ethm

Furthermore, we remark here that $A_p \to \infty$ as $p \to 1^+$ and we shall demonstrate this by an example, following the proof of Theorem \ref{KPSthm1}.

\br\label{KPS-rem1}
 It would be interesting to obtain an analog of Theorem \ref{KPSthm1} in multidimensional case. For instance, if $f:\mathbb{C} \to \mathbb{R}^n$ ($n \ge 3$) is harmonic
in  $\mathbb{D}$, then is it true that
$$\int_{-1}^{1}|f(x)|^p dx \leq  A_p  \int_{0}^{2\pi}|f(e^{i\theta})|^p d\theta?
$$

\er


%
Moreover, for the special case $p=2$, we get two variant forms of Theorem \ref{KPSthm1} and they are as follow:

\bthm\label{KPSthm3}
Let $f(z)=h(z)+ \overline{g(z)} = \sum _{k=0}^{\infty}a_k z^k + \sum _{k=1}^{\infty}\overline{b_k z^k} $ be harmonic in $\ID$. If $f \in {\bf h}^2$, then the following sharp inequality holds:
$$
\int_{-1}^{1}|f(x e^{it})|^2 dx \leq \frac{1}{2} \int_{0}^{2\pi}|f(e^{i\theta})|^2 d\theta
+2\pi {\rm Re} \left(\sum_{k=1}^\infty a_k b_k e^{i 2kt}\right) ~\mbox{ for all }~ t \in \IR.
$$
\ethm

\bthm\label{KPSthm4}
If $f=h+\overline g \in {\bf h}^2$,  then the following sharp inequalities hold for all $t \in \IR$:
\be\label{KPSeq8}
\int_{-1}^{1}|f(x e^{it})|^2 dx \leq \frac{1}{2} \int_{0}^{2\pi}|f(e^{i\theta})|^2 d\theta+||hg||_1
\leq  \int_{0}^{2\pi}|f(e^{i\theta})|^2 d\theta ~\mbox{ for all }~ t \in \IR.
\ee
\ethm

Proofs of all these results are presented in Section \ref{KPSSec2}.

\section{Proofs of Main Theorems and related results}\label{KPSSec2}
A positive real-valued function $u$ is called log-subharmonic, if $\log u$ is subharmonic. In order to prove our Theorems, we need the following classical inequality of Lozinski \cite{Lozinski} which is popularly known as
Fej\'er-Riesz-Lozinski  inequality, and also an inequality of Kalaj from \cite{Kalaj2017}.

\begin{Lem}\label{Lem_Lozin}
{\rm \cite{Lozinski}}
Suppose that $\Phi$ is a log-subharmonic function from $\ID$ to $\IR$, such that $\ds \int_{0}^{2\pi} \Phi^p(re^{i\theta})\, d\theta $ are uniformly bounded with respect to r for some $p>0$. Then the following sharp inequality holds:
\be\label{KPSeq4}
\int_{-1}^{1}\Phi^p(x e^{i t}) dx \leq \frac{1}{2} \int_{0}^{2\pi}\Phi^p(e^{i\theta}) d\theta ~~\mbox{ for all }~~ t \in \IR.
\ee
If equality is attained for some $t \in \IR$, then $\Phi \equiv 0$ in $\ID$. The constant $1/2$ is best possible.
\end{Lem}

As an application to Lemma \Ref{Lem_Lozin}, we deduce the following result.

\begin{lem}\label{KPSlem1}
Let $\varphi$ and $\psi$ be a pair of two analytic functions defined on $\ID$ such that $\varphi$ and $\psi \in H^p$ for  some $p > 1$. Then
\begin{equation} \label{Losin}
\int_{-1}^{1}(|\varphi(x e^{i t})|+|\psi(x e^{i t})|)^p dx \leq \frac{1}{2}\int_{0}^{2\pi}(|\varphi(e^{i\theta})|+|\psi(e^{i\theta})|)^p d\theta ~~\mbox{ for all }~~ t \in \IR.
\end{equation} The constant $1/2$ is sharp. \rm
\end{lem}
\bpf
In order to prove the result, it is enough to show that $\log (|\varphi(z)| + |\psi(z)|)$ is subharmonic in $\ID$. Then, we can apply Lemma \Ref{Lem_Lozin} and obtain the desired conclusion.
It is well known that $\log(|A(z)|^2 + |B(z)|^2)$ is subharmonic in $\ID$ provided $A(z)$ and $B(z)$ are analytic in $\ID$. Without loss of generality we can suppose that
both the functions $\varphi$ and $\psi$ are nonzero at each point of the unit disk. Then, there exist two non-vanishing analytic functions $A(z)$ and $B(z)$ in $\ID$ such that
$A^2(z) = \varphi(z)$ and $B^2(z)= \psi(z)$, which clearly implies that $\log (|\varphi(z)| + |\psi(z)|)$ is subharmonic in $\ID$.\\

Suppose that $\varphi (z)$ and $\psi(z)$ have zero(s) inside $\ID$. Then the zero(s) can be accumulated in the Blaschke product
so that
$$ \varphi (z) = B_1(z) A(z) ~~\mbox{ and }~~ \psi (z) = B_2(z) B(z),
$$
where $B_1(z)$, $B_2(z)$ are Blaschke products, $A(z)$ and $B(z)$ are non-vanishing analytic functions in $\ID$. Then, we deduce that
\beqq
\int_{-1}^{1}(|\varphi(x e^{i t})|+|\psi(x e^{i t})|)^p dx &=& \int_{-1}^{1}(|B_1(x e^{i t}) A(x e^{i t})|+|B_2(x e^{i t}) B(x e^{i t})|)^p dx \\
& \leq & \int_{-1}^{1}(|A(x e^{i t})|+|B(x e^{i t})|)^p dx \\
& \leq & \frac{1}{2} \int_{0}^{2\pi}(|A(e^{i\theta})|+|B(e^{i\theta})|)^p dx ~~ \mbox{(by Lemma \Ref{Lem_Lozin})}\\
& = & \frac{1}{2}\int_{0}^{2\pi}(|\varphi(e^{i\theta})|+|\psi(e^{i\theta})|)^p d\theta.
\eeqq
This completes the proof.
\epf

\begin{Thm}\label{thm_Kalaj_ana_har}
Let $1 < p < \infty$. Assume that $f = h +\overline{g} \in {\bf h}^p$ with ${\rm Re}(h(0)g(0))=0$. Then we have the following sharp inequality
\be\label{KPSeq5}
\int_{0}^{2\pi}(|h(e^{i \theta})|^2 +|g(e^{i \theta})|^2)^{p/2} d\theta \leq \frac{1}{\left(1 - |\cos \frac{\pi}{p}|\right)^{p/2}} \int_{0}^{2\pi}|f(e^{i\theta})|^p d\theta .
\ee
\end{Thm}

Theorem \Ref{thm_Kalaj_ana_har} has been stated and proved in \cite{Kalaj2017} and the proof uses plurisubharmonic function method initiated by  Hollenbeck and Verbitsky \cite{HV}. The following result due to Frazer is also useful in the proof of Theorem \ref{KPSthm1}.

\begin{Thm}\label{Frazer_subhar}{\rm \cite[Theorem 2.2]{Frazer}} If $U$ is subharmonic, positive, and continuous inside and on a circle $\Gamma$, and if $D_0, D_1,\ldots, D_{n-1}$ are $n$ diameters of $\Gamma$ such that the angles between the consecutive radii is $\pi/n$, then there exists a constant $A>0$ such that the following inequality holds:
$$\sum_{k=0}^{n-1}\int_{D_k}U^p(x)dx \leq  B_p  \int_{0}^{2\pi} U^p(e^{i\theta}) d\theta,
$$
where
$$B_p := \left\{
    \begin{array}{ll}
        \ds \frac{A n}{p-1}   & \mbox{if } 1 < p < 2 \\
        \ds \csc\left(\frac{\pi}{2n}\right)  & \mbox{if } p \geq 2.
    \end{array} \right.
$$
\end{Thm}

We remark here that Theorem \Ref{Frazer_subhar} is still valid, if we replace the condition ``positive'' by ``non-negative''.

\subsection{Proof of Theorem \ref{KPSthm1}}
Let $f=h+\overline{g} \in {\bf h}^p$ for some $p \in (1, 2]$. Without loss of generality, we can assume that $g(0)=0$. Then
 we have
\beqq
\int_{-1}^{1}|f(x e^{i t})|^p dx &\leq & \int_{-1}^{1}(|h(x e^{i t})|+|g(x e^{i t})|)^p dx \\
&\leq &
\frac{1}{2}\int_{0}^{2\pi}(|h(e^{i\theta})|+|g(e^{i\theta})|)^p d\theta \quad \mbox{(by Lemma \ref{KPSlem1})}\\
&\leq & \frac{2^{p/2}}{2}\int_{0}^{2\pi}(|h(e^{i\theta})|^2+|g(e^{i\theta})|^2)^{p/2} d\theta\\
&\leq & \frac{2^{p/2}}{2(1-|\cos(\pi/p)|)^{p/2}}\int_{0}^{2\pi}|f(e^{i\theta})|^p d\theta
\quad \mbox{(by Theorem \Ref{thm_Kalaj_ana_har})}.
\eeqq
It is a simple exercise to see that
$$ \frac{2^{p/2}}{2(1-|\cos(\pi/p)|)^{p/2}} = \frac{1}{2}\sec^{p} \left(\frac{\pi}{2 p}\right)  ~~~\mbox{ for } 1 < p \leq 2.
$$


Now, let us consider the case $p \geq 2$. Let $f \in {\bf h}^p$ for some $p >2$. Set $f_{\rho}(z) = f(\rho z)$. It is easy to see that $|f_{\rho}(z)|$ is subharmonic, non-negative,
and continuous inside and on the circle $|z| =1$ for $0 < \rho < 1$. Therefore, we can apply Theorem \Ref{Frazer_subhar} with $n=1$ and we get
\beqq
\int_{-\rho}^{\rho}|f(x e^{i t})|^p dx &=& \int_{-1}^{1}|f_{\rho}(x e^{i t})|^p dx \\
&\leq&   \int_{0}^{2\pi} |f_{\rho}(e^{i\theta})|^p d\theta \quad \mbox{(by Theorem \Ref{Frazer_subhar})}\\
&\leq&   \int_{0}^{2\pi} |f(e^{i\theta})|^p d\theta \quad \mbox{(by the monotonicity of $M_p(r, f)$)}.
\eeqq
Since the above inequality holds for all $\rho$ such that $0 < \rho < 1$, the desired conclusion follows.

Now, let us  prove the sharpness of the result for $p \in (1, 2]$ by an example.
By $K(p)$, we mean the optimal constant in the inequality
$$
\int_{-1}^{1}|f(x)|^p dx \leq  K(p)\int_{0}^{2\pi}|f(e^{i\theta})|^p d\theta, \quad p \ge 1.
$$
Let us show that
$$
K(p) \ge \frac{1}{2}\cos^{-p} \left(\frac{\pi}{2 p}\right), \quad p \ge 1.
$$
In particular $K(1)=+\infty$.
For $0 < r < 1$, consider the function
$$
f_r(z)={\rm Re} \left(\frac{1+r z}{1-r z}\right)^{1/p}, ~z\in\ID,
$$
where $r$ will be chosen as close to $1$. At first, it is easy to compute that
$$\int_{-1}^{1}|f_r(x)|^p dx=\frac{2}{r} \log \left( \frac{1+r}{1-r}\right)-2 =\frac{4}{r}{\rm arctanh}\, r -2.
$$
On the other hand
\beq\label{KPSeq6}
 \int_{0}^{2\pi}|f_r(e^{i\theta})|^p d\theta&=&\int_{0}^{2\pi}
 \left|\cos^p \left(\frac{1}{p}\arg\frac{1+re^{i\theta}}{1-re^{i\theta}} \right)\right|\,
 \left|\frac{1+re^{i\theta}}{1-re^{i\theta}}\right| d\theta \nonumber \\
 &=&\int_{0}^{2\pi}\left|\cos^p \left(\frac{1}{p}\arctan\frac{2 r\sin \theta}{1-r^2} \right)\right|
 \,\left|\frac{1+re^{i\theta}}{1-re^{i\theta}}\right| d\theta .
\eeq
We are interested in estimating the last integral when $r$ is close to $1$ and for that,
we need to consider the behaviour of the integrand when $r$ is close to $1$.
Given any arbitrarily small $\epsilon >0$ and large $M>0$, we can find a $\rho \in (0, 1)$ such that
$$\left|\frac{2r \sin \theta}{1-r^2} \right| > M ~\mbox{ for all }~ \theta \in E ~\mbox{ and }~ r \in (\rho, 1),
$$
where $E = [\epsilon/4, \pi-\epsilon/4] \bigcup [\pi + \epsilon/4, 2\pi-\epsilon/4]$.
It is a routine matter to see that
\beqq
I &=& \int_{0}^{2\pi}\frac{d\theta}{|1-re^{i\theta}|}\\
&=& 4\int_{0}^{\pi/2}\frac{d\theta}{\sqrt{(1+r)^2 -4r \sin^2\theta}}\\
&=& \frac{4}{1+r}\int_{0}^{\pi/2}\frac{d\theta}{\sqrt{1  - k^2\sin^2\theta}}, \quad k^2 = \frac{4r}{(1+r)^2},\\
&=& \frac{4{\K}(k)}{1+r},
\eeqq
where
$${\K}(k):= \frac{\pi}{2} \int_{0}^{\pi/2}\frac{d\theta}{\sqrt{1  - k^2\sin^2\theta}}
\quad (|k|<1)
$$
which is known as the complete elliptic integral of the first kind.
Therefore, for $r$ close to $1$, the integral in \eqref{KPSeq6} behaves like
$$2 \cos^{p} \left(\frac{\pi}{2 p}\right) \int_{0}^{2\pi}\frac{d\theta}{|1-re^{i\theta}|} =
8\cos^{p} \left(\frac{\pi}{2 p}\right)\left [ \frac{1}{1+r} {\K}\left(\frac{2\sqrt{r}}{1 + r}\right)\right ].
$$
A simple computation gives
$$
\ds \lim_{r \to 1}\left (\frac{\ds \int_{-1}^{1}|f_r(x)|^p dx}{\ds \int_{0}^{2\pi}|f_r(e^{i\theta})|^p d\theta}
\right )=\frac{1}{2}\cos^{-p} \left(\frac{\pi}{2 p}\right),
$$
which proves that the bounds in Theorem \ref{KPSthm1} are sharp for $1<p\leq 2$.
The proof is complete. \hfill{$\Box$}

 The example demonstrated in the proof of Theorem \ref{KPSthm1} suggests the following.

 \bcon
 If $f \in {\bf h}^p$ for $p>2$, then the following sharp inequality holds:
 $$
 \int_{-1}^{1}|f(x e^{i t})|^p dx \leq  \frac{1}{2}\cos^{-p} \left(\frac{\pi}{2 p}\right) \int_{0}^{2\pi}|f(e^{i\theta})|^p d\theta ~~\mbox{ for all }~~ t \in \IR.
 $$
\econ

Next, we present the variant forms of Riesz-Fej\'er theorems for harmonic functions for the special case $p=2$.
\subsection{Proof of Theorem \ref{KPSthm3}}
It is enough to prove the case $t=0$ as $F_t(z) = f(z e^{it}) \in {\bf h}^2$ for every fixed $t \in \IR$. Let
$$\Phi(z)=h(z)+\overline{g(\overline z)} = \sum _{k=0}^{\infty}a_k z^k + \sum _{k=1}^{\infty}\overline{b_k} z^k.$$
Then, it is easy to see that
\beqq
\int_{-1}^{1}|f(x)|^2 dx &=& \int_{-1}^{1}|\Phi(x)|^2 dx \\
&\leq & \frac{1}{2} \int_{0}^{2\pi}|\Phi(e^{i\theta})|^2 d\theta\\
& = & \frac{1}{2} \int_0^{2\pi} (|h(e^{i\theta})|^2+|g(e^{i\theta})|^2) d\theta+ {\rm Re} \int_0^{2\pi} h(e^{i\theta}) g(e^{-i\theta})d\theta \\
&=& \frac{1}{2} \int_0^{2\pi} (|h(e^{i\theta})|^2+|g(e^{i\theta})|^2) d\theta+2\pi {\rm Re} \sum_{k=0}^\infty a_k b_k \\
&=& \frac{1}{2} \int_{0}^{2\pi}|f(e^{i\theta})|^2 d\theta+2\pi {\rm Re} \sum_{k=1}^\infty a_k b_k \quad \mbox{(since $b_0=0$)}.
\eeqq
This completes the proof. \hfill{$\Box$}
\vspace{8pt}

\subsection{Proof of Theorem \ref{KPSthm4}}
Let $f = h + \overline{g} \in {\bf h}^2$. Without loss of the generality we may suppose that $g(0)=0$. Furthermore, it is enough to consider the case $t=0$ as the other cases follow in a similar way. Accordingly,
\beqq
\int_{-1}^{1}|f(x)|^2 dx &=& \int_{-1}^{1}|h(x)|^2 dx+\int_{-1}^{1}|g(x)|^2 dx+2{\rm Re} \int_{-1}^{1}h(x)g(x) dx\\
& \leq & \int_{-1}^{1}|h(x)|^2 dx+\int_{-1}^{1}|g(x)|^2 dx+2 \int_{-1}^{1}|h(x)g(x)| dx.
\eeqq
As $h(z)g(z)$ is analytic in $\mathbb{D}$, from the Cauchy-Schwarz inequality, we observe that $hg \in H^1$, and hence from
Riesz-Fej\'{e}r inequality we get
\beqq
\int_{-1}^{1}|f(x)|^2 dx & \leq & \frac{1}{2}\int_{0}^{2\pi}|h(e^{i\theta})|^2 d\theta+ \frac{1}{2} \int_{0}^{2\pi}|g(e^{i\theta})|^2 d\theta+ \int_{0}^{2\pi}|h(e^{i\theta})g(e^{i\theta})| d\theta\\
& = & \frac{1}{2}\int_{0}^{2\pi}|f(e^{i\theta})|^2 d\theta + ||hg||_1.
\eeqq
Finally, the right hand side inequality in \eqref{KPSeq8} follows from the following inequality
$$
||hg||_1 \leq ||h||_2||g||_2 \leq \frac{1}{2}(||h||_2^2+||g||_2^2)= \frac{1}{2}\int_{0}^{2\pi}|f(e^{i\theta})|^2 d\theta.
$$

Now let us show that both the inequalities in Theorem \ref{KPSthm4} are sharp and cannot be improved.
For this, we consider
$$f(z) = h(z)+\overline {g(z)} = \frac{\sqrt{\varphi'(z)}}{2} + \overline{\frac{\sqrt{\varphi'(z)}}{2}},
$$
where $\varphi$ maps the disk $\mathbb{D}$ onto the rectangle with the vertices $\pm 1 \pm i \varepsilon$
such that $[-1,1]$ maps onto itself. Let us remark that $\varphi'(x)>0$ and $\varphi'(0) \to 0$ as $\varepsilon \to 0$.
From the basic calculus, it is clear that
$$
\int_{-1}^{1}|f(x)|^2 dx=\int_{-1}^{1}|\varphi'(x)| dx=2.
$$

On the other hand
\beqq
\int_{0}^{2\pi}|f (e^{i\theta})|^2 d\theta &=& \int_{0}^{2\pi}|{\rm Re} \sqrt{\varphi' (e^{i\theta})}|^2 d\theta \\
&=& \frac{1}{2}\int_{0}^{2\pi}| \varphi' (e^{i\theta})|d\theta+\pi {\rm Re}( \varphi'(0)) \\
&=& \frac{4+4\varepsilon}{2}+\pi {\rm Re}( \varphi'(0)) \to 2 \mbox{ as }\varepsilon \to 0.
\eeqq

Furthermore, as $\varepsilon \to 0$, we get
$$ ||hg||_1 = \int_{0}^{2\pi}|h(e^{i\theta})g(e^{i\theta})| d\theta = \frac{1}{4} \int_{0}^{2\pi}|\varphi'(e^{i\theta})| d\theta
= \frac{1}{2} \int_{0}^{2\pi}|f (e^{i\theta})|^2 d\theta.
$$
So, from here we see that the constant on the right hand side of the inequality cannot be less than 1 and both the inequalities are sharp. The proof
is complete. \hfill{$\Box$}
%
\subsection*{Acknowledgements}
The research of the first author was funded by the subsidy allocated to Kazan Federal University for the state assignment in the sphere of scientific activities, project No. 1.12878.2018/12.1, and the work of the second author is supported by Mathematical Research Impact Centric Support
(MATRICS) grant, File No.: MTR/2017/000367, by the Science and Engineering Research Board (SERB), Department of Science and Technology (DST), Government of India. The third author initiated this work when he was a NBHM Postdoctoral Fellow at the Indian Statistical Institute Chennai Centre.

\end{document}